\documentclass[10pt]{amsart}
\usepackage{amsmath}
\usepackage{amssymb}
\usepackage{amsthm}
\usepackage{amsfonts}

\def\ad{\textup{ad} }
\def\Ad{\textup{Ad } }
\def\H{{\mathcal H}}
\def\C{\mathbb{C}}

\def\dim{\textup{dim }}
\def\b{\mathfrak{b}}

\newcommand{\g}[1]{{\mathfrak{g}}_{#1}}
\newcommand{\n}[1]{{\mathfrak{n}}_{#1}}

\theoremstyle{plain}
\newtheorem{theorem}{Theorem}
\newtheorem{lemma}[theorem]{Lemma}
\newtheorem{proposition}[theorem]{Proposition}
\newtheorem{corollary}[theorem]{Corollary}

\newtheorem*{theorem*}{Theorem}
\numberwithin{theorem}{section}

\theoremstyle{definition}
\newtheorem{definition}[theorem]{Definition}

\numberwithin{equation}{section}

\title{Paving Hessenberg varieties by affines}
\author{Julianna S. Tymoczko}
\address{Department of Mathematics, University of Iowa, Iowa City, IA 52242}
\email{tymoczko@math.uiowa.edu}

\begin{document}

\begin{abstract}
Regular nilpotent Hessenberg varieties form a family of subvarieties
of the flag variety arising in the study of quantum cohomology,
geometric representation theory, and numerical analysis.
In this paper we construct a paving by affines
of regular nilpotent Hessenberg varieties for all classical types, generalizing 
results of de Concini-Lusztig-Procesi and Kostant.
This paving is in fact the intersection of a particular Bruhat
decomposition with the Hessenberg variety.  The nonempty cells of
the paving and their dimensions are identified by 
combinatorial conditions on roots.
We use the paving to prove these Hessenberg varieties have no 
odd-dimensional homology.  
\end{abstract}

\maketitle

\section{Introduction}

This paper studies the topology of regular nilpotent Hessenberg varieties, a family of subvarieties of the flag variety introduced in \cite{dMPS} that arise naturally in contexts as diverse as numerical analysis, number theory, and representation theory.  Their geometry encodes deep algebraic and combinatorial properties, including the quantum cohomology of the flag variety \cite{Ko}.  We prove that regular nilpotent Hessenberg varieties in classical Lie types have a {\em paving by affines}, a cell decomposition like CW-decompositions but with weaker closure relations.  This paving permits us to describe the varieties' cohomology, for instance to show that it vanishes in odd dimensions.  Moreover, this paving can be realized as the intersection of the Hessenberg variety with a particular Bruhat decomposition, so the dimensions of the nonempty cells are characterized by combinatorial conditions.

Let $G$ be a complex linear algebraic group of classical type, $B$ a fixed Borel subgroup, and $\g{}$ and $\b$ their Lie algebras.  A Hessenberg space $H$ is a linear subspace of $\g{}$ that contains $\b$ and that is closed under Lie bracket with $\b$, namely $[H,\b]$ is contained in $H$.  Fix an element $X$ in $\g{}$ and a Hessenberg space $H$. The Hessenberg variety $\H(X,H)$ is the subvariety of the flag variety $G/B$ consisting of $gB/B$ satisfying $g^{-1}Xg \in H$, or equivalently $\Ad g^{-1}(X) \in H$.

An important special case is when $G = GL_n(\C)$, $B$ consists of the upper-triangular invertible matrices, $\g{}$ is the set of all $n \times n$ matrices, and $\b \subseteq \g{}$ the subset of all upper-triangular matrices.  In this case, a Hessenberg space $H$ is equivalent to a nondecreasing function $h: \{1,2,\ldots, n\} \rightarrow \{1,2,\ldots,n\}$ satisfying $h(i) \geq i$ for all $i$, by the rule that $H$ is a subspace of $\g{}$ whose matrices vanish in positions $(i,j)$ whenever $i > h(j)$.  The flags in $GL_n(\C)/B$ can be written as nested subspaces $V_1 \subseteq V_2 \subseteq \cdots \subseteq V_{n-1} \subseteq \C^n$ where each $V_i$ is $i$-dimensional.  The Hessenberg variety $\H(X,H)$ is the collection of flags for which $X V_i \subseteq V_{h(i)}$ for each $i$.

Let $N$ be a regular nilpotent element of $\g{}$, namely let $N$ be in the dense adjoint orbit within the nilpotent elements in $\g{}$.  When $G=GL_n(\C)$, the regular nilpotent elements are those which consist of a single Jordan block.  In this paper we prove the existence of a paving by affines for regular nilpotent Hessenberg varieties.  Pavings, defined formally in Section \ref{pavings}, are like CW-decompositions but with weaker closure relations.  The paving in this paper can be described explicitly.

\begin{theorem*}
Fix a Hessenberg space $H$ with respect to the Borel $\mathfrak{b}$
and fix a regular nilpotent element $N$ in $\mathfrak{b}$.  
The Bruhat decomposition $BwB/B$ of the flag variety intersects
the Hessenberg variety $\H(N,H)$ in a paving by affine cells.
The cell $P_w = \H(N,H) \cap BwB/B$ corresponding to $w$ 
is nonempty if and only if
$\Ad w^{-1} (E_{\alpha_j})$ is in $H$ for each simple root vector
$E_{\alpha_j}$.  If $P_w$ is nonempty its dimension is given by
$\dim (\mathfrak{b} \cap \Ad w (\mathfrak{b}^- \cap H)) - \textup{rank}(G)$, where
$\mathfrak{b}^-$ denotes the opposite Borel.
\end{theorem*}

Theorem \ref{main theorem} gives the complete statement of the paving result, including
conditions on roots which determine when $P_w$ is nonempty and, if so, its dimension.  The arbitrary regular nilpotent Hessenberg variety $\H(N,H)$ is also paved by affines (Corollary \ref{generic paving}), since it is homeomorphic to an $\H(N', H)$ satisfying the conditions of the Theorem. Using this paving, Corollary \ref{zero odd} proves regular nilpotent Hessenberg varieties have no odd-dimensional cohomology.

This result generalizes from type $A_n$ some of the work of \cite{T}, where the reader may find explicit examples for $GL_n(\C)$.  It also extends results of \cite{dCLP} beyond the Springer fiber, namely when $H = \b$ and $X$ is arbitrary.  (The cohomology of the Springer fiber carries a natural Weyl group action which geometrically constructs all the irreducible representations of the Weyl group; see \cite{S}, \cite{BM}, \cite{L}, and \cite{CG}, among others.)  The Theorem strengthens work of B.~Kostant \cite{Ko} for the Peterson variety, i.e. when $H$ is generated by $\b$ together with the root spaces corresponding to the negative simple roots.  In \cite{Ko}, Kostant intersected the Peterson variety with a different Bruhat decomposition that paved it by affine varieties rather than the affine cells $\mathbb{C}^k$ used here.  (Kostant's paving gave an open dense subvariety of the Peterson variety with coordinate ring isomorphic to the quantum cohomology ring of the flag variety.)  The topology of regular nilpotent Hessenberg varieties also gives information about algebraic invariants of ad-nilpotent ideals in a Borel subalgebra, which are closely related to Hessenberg spaces \cite{ST}.
This paper parallels the results for regular semisimple Hessenberg varieties in \cite{dMPS} using a different approach. We understand D.~Peterson has uncirculated results overlapping these \cite{C}, \cite[Theorem 3]{BC}; we infer that the methods used here are substantively different.  

The second section has background information, including the definition of pavings, Bruhat decompositions, and the decomposition of the nilradical of $\g{}$ into subspaces called {\em rows}.  The third section identifies the restriction and projection of the map $\ad N$ to the individual rows.  The fourth uses these results to show that the cells of this Bruhat decomposition intersect the Hessenberg variety $\H(N,H)$ in an iterated tower of affine fiber bundles.  

This paper was part of the author's doctoral dissertation.  The author gratefully acknowledges the comments from and helpful discussions with Mark Goresky, David Kazhdan, Arun Ram, Konstanze Rietsch, and Eric Sommers, and especially those of the author's advisor, Robert MacPherson.

\section{Background and definitions}

This section contains background and definitions needed from the
literature for the rest
of the paper.  The first subsection recalls the necessary results
about pavings.
The second subsection reviews the Bruhat decomposition 
as well as parameterizations of each Schubert cell.
The third subsection defines a partition of the positive roots
into {\em rows}.  These rows span subspaces of $\b$ which are
abelian or Heisenberg, simplifying later computations.

\subsection{Pavings} \label{pavings}

Pavings are common decompositions of algebraic varieties.

\begin{definition}
A paving of an algebraic variety $X$ is an ordered partition into disjoint
$X_0$, $X_1$, $X_2$, $\ldots$ so that each finite union
$\bigcup_{i=0}^j X_i$ is Zariski-closed in $X$.
\end{definition}

The $X_i$ are the cells of the paving.  Note that
pavings have weaker closure relations than CW-decompositions since
the boundary of a cell is not required to be contained in cells of
smaller dimension.

\begin{definition}
A paving by affines of $X$ is a paving so that each $X_i$ is 
homeomorphic to affine space.
\end{definition}

The following is the main reason we use pavings \cite[19.1.11]{F}.

\begin{lemma} \label{odd-dim}
Let $X=\bigcup X_i$ be a paving by a finite number of affines with
each $X_i$ homeomorphic to $\C^{d_i}$.  The cohomology groups
of $X$ are given by $H^{2k}(X) = \bigoplus_{i:d_i = k} \mathbb{Z}$.
\end{lemma}

\subsection{Bruhat decompositions}

Fix a Borel subgroup $B$ and a maximal torus $T$ in $B$, 
and let $W$ denote the Weyl group of $G$, namely
the quotient $N(T)/T$ of the normalizer of $T$.
The subgroup $B$ determines a decomposition of the flag variety
$G/B$ into cosets $BwB/B$ indexed by the elements $w$ of the Weyl group $W$.  

In fact, this is a classic paving by affines.  Recall that 
the length of the element $w$ is the minimal
number of simple transpositions $s_1$, $\ldots$, $s_n$
required to write $w = s_{i_1} \cdots s_{i_{\ell(w)}}$.  The next lemma is proven in \cite{Ch}, among others.

\begin{lemma}
The cells $BwB/B$ of the Bruhat decomposition form a paving by affines
when ordered in any way subordinate to the partial order determined by
the length of $w$.
\end{lemma}

We use an explicit description of the affine cells of this
paving.
Let $\Phi$ denote the roots of $\g{}$ and $\Phi^+$ the roots corresponding
to $\b$.  Recall the partial order on $\Phi$ given by $\alpha > \beta$
if and only if $\alpha - \beta$ is a sum of positive roots.  Write
$\g{\alpha}$ for the root space corresponding to $\alpha$,
$U$ for the maximal unipotent subgroup of $B$, $U^-$ for
its opposite subgroup, and $\n{}$ for the Lie algebra of $U$.

\begin{lemma}
Fix $w$ in $W$.  The following are homeomorphic:
\begin{enumerate}
\item the Schubert cell $BwB/B$;
\item the subgroup $U_w = \{u \in U: w^{-1}uw \in U^{-}\}$;
\item the Lie subalgebra $\n{w} = \bigoplus_{\alpha \in \Phi^+: 
w^{-1} \alpha < 0} \g{\alpha}$.
\end{enumerate}
\end{lemma}

\begin{proof}
The subgroup $U_w$ forms a set of coset representatives for $BwB/B$
and is a product of root subgroups
$U_w = \prod_{w^{-1} \alpha < 0} U_{\alpha}$ for any fixed order of the roots
\cite[Theorems 28.3 and 28.4, Proposition 28.1]{H}.  
Since $\n{w}$ is nilpotent its image under the exponential map is $\exp{\n{w}}=U_w$, 
as in \cite[page 50]{K}.
\end{proof}

Let $\Phi_w = \{\alpha \in \Phi^+: w^{-1} \alpha < 0\}$ 
be the set of roots indexing $U_w$ and $\n{w}$.

\subsection{Rows}

This subsection describes a partition of positive roots
into {\em rows} which facilitates inductive proofs
because each row generates an abelian or Heisenberg subalgebra of $\g{}$.  
The subsection also includes a table enumerating the roots in each row
of classical type.  A version of 
this decomposition is used elsewhere, e.g. \cite{Ste}.  

We are motivated by $GL_n(\C)$, where the unipotent group
$U$ can be taken to be upper-triangular matrices with ones along the
diagonal.  In this case, the $i^{th}$ row corresponds to the subgroup
of $U$ with nonzero entries only along the $i^{th}$ row and the diagonal.
Direct computation shows that this subgroup is abelian and that the
product of the rows is $U$.

Write $\alpha \geq \beta$ to indicate  
either $\alpha > \beta$ or $\alpha = \beta$.
The roots of the $i^{th}$ row are
\[\Phi_i = \{ \alpha \in \Phi^+: \alpha \geq \alpha_i, \alpha \not > 
\alpha_j \textup{ for each } j=1, \ldots i-1 \}.\]

The following labelling of the simple roots in classical types

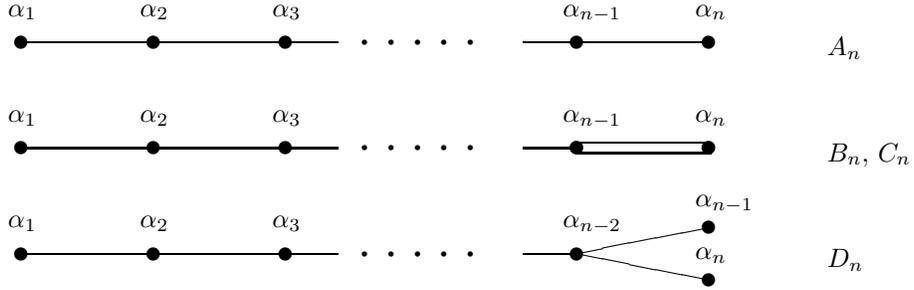
\begin{picture}(330,120)(100,-50)
\multiput(95,45)(50,0){3}{\circle*{5}}
\multiput(95,5)(50,0){3}{\circle*{5}}
\multiput(95,-35)(50,0){3}{\circle*{5}}
\multiput(225,45)(10,0){5}{\circle*{2}}
\multiput(305,45)(50,0){2}{\circle*{5}}
\put(400,40){$A_n$}
\multiput(225,5)(10,0){5}{\circle*{2}}
\multiput(305,5)(50,0){2}{\circle*{5}}
\multiput(225,-35)(10,0){5}{\circle*{2}}
\multiput(305,5)(50,0){2}{\circle*{5}}
\put(400,0){$B_n$, $C_n$}
\put(400,-40){$D_n$}
\multiput(95,45)(50,0){2}{\line(1,0){50}}
\multiput(95,5)(50,0){2}{\line(1,0){50}}
\multiput(95,-35)(50,0){2}{\line(1,0){50}}
\multiput(195,45)(0,-40){3}{\line(1,0){20}}
\put(305,-35){\circle*{5}}
\multiput(355,-25)(0,-20){2}{\circle*{5}}
\put(305,45){\line(1,0){50}}
\put(305,7){\line(1,0){51}}
\put(305,3){\line(1,0){51}}
\put(305,-35){\line(5,1){50}}
\put(305,-35){\line(5,-1){50}}
\multiput(305,45)(0,-40){3}{\line(-1,0){20}}
\multiput(90,55)(0,-40){3}{$\alpha_1$}
\multiput(140,55)(0,-40){3}{$\alpha_2$}
\multiput(190,55)(0,-40){3}{$\alpha_3$}
\multiput(300,55)(0,-40){2}{$\alpha_{n-1}$}
\multiput(350,55)(0,-40){2}{$\alpha_{n}$}
\put(300,-25){$\alpha_{n-2}$}
\put(350,-17){$\alpha_{n-1}$}
\put(350,-37){$\alpha_{n}$}
\end{picture}

gives the partition into rows of Table \ref{classical rows},
which is used throughout this paper.

\begin{figure}[h]
\begin{tabular}{c|l}
Type & \hspace{35mm} Row $\Phi_i$ \\
\cline{1-2} $A_n$ & $\{\sum_{j=i}^k \alpha_j: i \leq k \leq n\}$ \\
$B_n$ &  $\left\{\sum_{j=i}^k \alpha_j: i \leq k \leq n\right\}
\bigcup \left\{\sum_{j=k}^n \alpha_j + 
\sum_{j=i}^n \alpha_j: i+1 \leq k \leq n\right\}$ \\
$C_n$ & $\left\{\sum_{j=i}^k \alpha_j: i \leq k \leq n\right\}
\bigcup \left\{\sum_{j=k}^{n-1} \alpha_j + 
\sum_{j=i}^n \alpha_j: i \leq k \leq n-1\right\}$ \\
$D_n$ & $
\left\{\sum_{j=i}^k \alpha_j: i \leq k \leq n-1 \right\}$ \\ 
& \hspace{5em} $\bigcup \left\{\sum_{j=i}^{n-2} \alpha_j + \alpha_n + 
\sum_{j=k}^{n-1} \alpha_j: i+1 \leq k \leq n\right\}$
\end{tabular} \caption{Rows in classical types} \label{classical rows}
\end{figure}

The rows generate subalgebras of the Lie algebra which we also call
rows.  We use $\n{i}$ to denote the subalgebra spanned by the root
spaces corresponding to the roots of $\Phi_i$, so
\[\n{i} = \bigoplus_{\alpha \in \Phi_i} \g{\alpha}.\]

Recall that $\g{}$ is an abelian Lie algebra if 
the derived algebra 
$[\g{},\g{}]$ is zero.  We call $\g{}$ a Heisenberg Lie algebra
if its lower central series $\g{} \supsetneq [\g{},\g{}] \supsetneq
[\g{},[\g{},\g{}]]=0$ vanishes after two steps,
if its derived algebra $[\g{},\g{}]$ is a one-dimensional subalgebra,
and if for all $X$ in $\g{}$ the map 
$\ad X$ surjects onto $[\g{},\g{}]$ unless $X$ is in $[\g{},\g{}]$.
The next proposition gives a family of examples of Heisenberg
algebras.

\begin{proposition}
In type $A_n$, $B_n$, or $D_n$ each row $\n{i}$ is abelian.  In type $C_n$ the row $\n{i}$ is Heisenberg when $i$ is not $n$ and is abelian when $i=n$.
\end{proposition}

\begin{proof}
Both proofs rely on the property that
\begin{equation} \label{adjoint on root spaces}
[\g{\alpha},\g{\beta}] = \left\{ \begin{array}{ll}
    0 & \textup{ if } \alpha + \beta \textup{ is not a root, and} \\
    \g{\alpha+\beta} & \textup{ if } \alpha + \beta \textup{ is a root.}
    \end{array}\right.\end{equation}
By inspection of Figure \ref{classical rows} we see that
for no choice of $\alpha$, $\beta$ in 
$\Phi_i$ in types $A_n$, $B_n$, and $D_n$
is the sum $\alpha+\beta$ a root.  This implies that $\n{i}$ is abelian.  (The argument applies to $\Phi_n = \{\alpha_n\}$ in type $C_n$ as well.)

In type $C_n$ let $\gamma_i$ denote the root 
$\sum_{j=i}^{n-1} 2 \alpha_j + \alpha_n$.  If $\alpha$ is any root
in $\Phi_i$ other than $\gamma_i$ then
the difference $\gamma_i - \alpha$ is a root in $\Phi_i$.
This is the Heisenberg property.  Since $\gamma_i+\gamma_i$ is not a root,
the subalgebra $\n{i}$ is Heisenberg.
\end{proof}

We define $U_i$ to be the subgroup associated to $\n{i}$.  The group 
$U_i$ can be characterized either as the product $\prod_{\alpha \in \Phi_i} 
U_{\alpha}$ or as the exponential $\exp(\n{i})$.  

\begin{proposition}
The unipotent group $U$ factors as the product $U= U_1 U_2 \cdots U_n$.
\end{proposition}

\begin{proof}
$U = \prod_{\alpha \in \Phi^+} U_{\alpha}$ for any fixed ordering
of the positive roots \cite[Proposition 28.1]{H}.  
The rows are abelian possibly up to a root subgroup, which can be 
ordered last.
\end{proof}

The intersection of a Schubert cell with a row is
$(U_i \cap U_w) w B/B$ or equivalently $\exp(\n{i} \cap \n{w}) w B/B$.
 
\section{Adjoint actions on rows}

We now begin our study of adjoint actions on regular nilpotent elements.
Fix a regular nilpotent $N$ in $\n{}$ and a Hessenberg space $H$ with
respect to $\mathfrak{b}$.
For each group element $g$, 
we will choose an appropriate $u$ in $U$ and 
reduce the problem of determining if $\Ad g^{-1}(N)$ is in $H$ 
to the question of whether a summand of $\Ad u^{-1}(N)$ is in a fixed 
subspace of $\n{i}$.  This will rely on the key fact that the adjoint
representation of a row is ``almost'' linear,  made precise in
Proposition \ref{near-linearity}.

Fix a root vector $E_{\alpha}$ to generate the root space $\g{\alpha}$
and define $m_{\alpha, \beta}$ by 
$[E_{\alpha}, E_{\beta}] = m_{\alpha, \beta}
E_{\alpha + \beta}$.  By Equation \eqref{adjoint on root spaces}, the
coefficient $m_{\alpha, \beta}$ is nonzero if and only if 
$\alpha + \beta$ is a root.
The set $\{E_{\alpha}: \alpha \in \Phi^+\}$ form a basis for the 
Lie algebra $\n{}$.  We refer to the expansion $Y = \sum y_{\alpha}
E_{\alpha}$ as the basis vector expansion of $Y$.

Let $\rho_i: \n{} \longrightarrow \n{i}$ be the vector space projection
determined by this basis of root vectors.  For
each $N$ in $\n{}$, 
define the map $\theta_i(N): \n{} \longrightarrow \n{i}$ by the equation
$\theta_i(N)(X) = \rho_i \Ad \exp X(N)$.

\begin{proposition} \label{near-linearity}
Fix $N$ in $\n{}$ and $X$ in $\n{j}$.
\begin{enumerate}
\item If $i<j$ then
$\theta_i(N)(X) = \rho_i\left(N - \ad(N)(X) + \frac{\ad^2(X)(N)}{2}\right)$;
\item if $i=j$ then
\[\theta_i(N)(X) = \left\{ \begin{array}{ll}
    \rho_i \left(N - \ad(N)(X)\right)
    & \textup{ in types } A_n, B_n, D_n, \textup{ and } \\

    \rho_i\left(N -\ad(N)(X) + \frac{\ad^2(X)(N)}{2}\right) 
    & \textup{ in type } C_n; \end{array} \right.\]
\item and if $i > j$ then $\theta_i(N)(X) = \rho_i N$.
\end{enumerate} 
Furthermore, when $i=j$ in type $C_n$ the image of $\rho_i \ad^2 (X)$ lies in
$\g{\gamma_i}$.
\end{proposition}

\begin{proof}
Recall that $\Ad(\exp X) = \exp(\ad X) = \sum_{n=0}^{\infty} \frac{(\ad X)^n}{n!}$ as in \cite[Proposition 1.93]{K}.  Write $X = \sum_{\alpha \in \Phi_j} x_{\alpha} E_{\alpha}$ and $Y = \sum_{\alpha \in \Phi^+} y_{\alpha} E_{\alpha}$ in terms of the basis.  The adjoint operator $\ad(X)(Y) = [X,Y]$ can be expanded as $\ad(X)(Y) = \sum_{\alpha+ \beta \in \Phi} m_{\alpha, \beta} x_{\alpha}y_{\beta} E_{\alpha+\beta}$ for the nonzero coefficients $m_{\alpha, \beta}$ by Equation \eqref{adjoint on root spaces}.

The rest of the proof follows from this relation.
Each element in the image of $\ad^3 X$ is a linear
combination of root
vectors $E_{\beta}$ with $\beta \geq 3 \alpha_j$, since $X$ is
in $\n{j}$.  By Table \ref{classical rows}
no such root $\beta$ exists, so $\theta_i(N)$ is a polynomial 
of degree at most two in the $x_{\alpha}$.  

Now let $i=j$. 
In types $A_n$, $B_n$, and $D_n$ there is no root in $\Phi_i$ 
greater than $2 \alpha_i$, so $\ad^2 X$ vanishes.  In these
types $\theta_i(N)$ is affine.  In type
$C_n$ there is a unique root greater than $2 \alpha_i$, namely the
 root $\gamma_i = \sum_{k=i}^{n-1} 2 \alpha_k + \alpha_n$.  
The image of $\rho_i \ad^2 X$ must thus be in $\g{\gamma_i}$.

Finally, choose $i>j$ and 
let $c_{\beta} E_{\beta}$ be a nonvanishing summand in the 
expansion of $\sum_{n=1}^{\infty} \ad^nX (N)$.
Then $\beta$ is the sum $\beta_1 + \cdots + \beta_k + \alpha$ for 
$\beta_1$, $\ldots$, $\beta_k$ in $\Phi_j$.  By definition $\beta$ 
is contained one of the rows $\Phi_1$, $\Phi_2$,
$\ldots$, or $\Phi_j$.  This means that $\theta_i(N)(X) = \rho_i N$.
\end{proof}

We include the next lemma for ease of reference.  It is a restatement
of known results.

\begin{lemma} \label{reg nilps}
The element $N=\sum_{\alpha \in \Phi^+} n_{\alpha} E_{\alpha}$ 
is a regular nilpotent element of $\n{}$ if and only if $n_{\alpha_i}$ is
nonzero for each simple root $\alpha_i$.  The set of regular nilpotent
elements of $\n{}$ is exactly the orbit $\Ad B(N)$ for each regular
$N$ in $\n{}$.
\end{lemma}

\begin{proof}
Use \cite[Lemma 4.1.4]{CM} for regular nilpotents.
The lemma says in this case that $\Ad B(N)$ is exactly the set of regular
nilpotents in $\n{}$, and that $\Ad B(N)$ is the set 
$\Ad T(N) + [\n{},\n{}]$ where $T$ is the maximal torus in $B$.
Since $N$ can be taken to be the sum of the simple root vectors
by \cite[Theorem 4.1.6]{CM}, the claim follows.
\end{proof}

The element $N$ in $\n{}$ defines a map
in the endomorphism ring $\textup{End}(\n{i})$ that is the 
linear part of $\theta_i(N)$ when $\theta_i(N)$ is affine.

\begin{definition}
For each $N$ in $\n{}$ the map $\psi_i(N)$ is defined as
the restriction and projection 
$\psi_i(N) = (\rho_i \circ \ad N) |_{\n{i}}$.  
\end{definition}

As before, let $\{E_{\alpha}: \alpha \in \Phi^+\}$ be a 
fixed basis of root vectors for $\g{}$ with structure
constants given by $[E_{\alpha}, E_{\beta}]=
m_{\alpha, \beta} E_{\alpha + \beta}$.  
Write $N$ in terms of this basis as 
$N = \sum_{\alpha \in \Phi^+} n_{\alpha} E_{\alpha}$.

The next lemma establishes properties of $\psi_i(N)$ with respect to this basis, where the linear map $\psi_i(N)$ is identified with its matrix.  The entries of this matrix are indexed by pairs of roots $(\alpha, \beta)$ in $\Phi_i \times \Phi_i$.  For instance, the entry at position $(\alpha, \beta)$ is the coefficient of $E_{\alpha}$ in $\psi_i(N)(E_{\beta})$.

\begin{lemma} \label{matrix entries}
Fix $N$ in $\n{}$.

\begin{enumerate}
\item The $(\alpha, \beta)$ position of
$\psi_i(N)$ has entry $m_{\alpha - \beta, \beta} n_{\alpha - \beta}$ 
if $\alpha - \beta$ is a positive root and zero otherwise.

\item If $X$ is in $\n{i-j}$ for some positive $j$ then 
$\psi_i(\Ad \exp X (N)) = \psi_i(N)$.
\end{enumerate}
\end{lemma}

\begin{proof}
The first part follows from the construction of the basis.

The second part follows from the first 
once we identify the coefficients of
$E_{\alpha - \beta}$ in $\Ad \exp X (N)$, for each pair 
of roots $\alpha$ and $\beta$ in $\Phi_i$.  If the
difference $\alpha - \beta$ is a root then it must be in a row
indexed by $k$, where $k$ is at least $i$.  
The coefficient of $E_{\alpha-\beta}$ in
$\Ad \exp X(N)$ is the same as that in
$\rho_k \Ad \exp X(N)$, which is $n_{\alpha - \beta}$ 
by Proposition \ref{near-linearity}.
\end{proof}

\begin{corollary}
Fix a regular nilpotent $N$ in $\n{}$ in types $A_n$, $B_n$, or $C_n$.
The map $\psi_i(N)$ is a regular nilpotent element of $\textup{End}(\n{i})$.
\end{corollary}

\begin{proof}
Order the basis $\{E_{\alpha}: \alpha \in \Phi_i\}$
by the height of $\alpha$ from highest to lowest.  By 
Table \ref{classical rows} 
this is a total order in which each
root differs by a simple root from the next.

Consider the matrix for $\psi_i(N)$ with respect to this basis.
The entries on and below the diagonal correspond to
differences $\alpha - \beta$ which are not positive.  The matrix
for $\psi_i(N)$ is zero in these positions by Lemma \ref{matrix entries}.
The $(\alpha, \beta)$ position is immediately 
above the diagonal if and only if $\alpha$ is immediately before 
$\beta$ in the height order.  In this case $\alpha - \beta$ is a
simple root $\alpha_j$ and the corresponding entry of
$\psi_i(N)$ is $m_{\alpha_j, \alpha -  \alpha_j} n_{\alpha_j}$.  
This is nonzero because
$N$ is regular nilpotent, by Lemma \ref{reg nilps}.
Since $\psi_i(N)$ is an upper-triangular matrix with nonzero entries
above the diagonal, it too is regular nilpotent, using Lemma 
\ref{reg nilps} for $\mathfrak{gl}_n$.
\end{proof}

The following lemma is similar to the previous and is necessary to
handle technical difficulties in type $D_n$, where
$\psi_i(N)$ is not a regular nilpotent operator.  
Another analogue of the previous lemma for type $D_n$ is given in 
Lemma \ref{main lemma D}.

\begin{lemma} \label{D_n coefficients}
Fix $X = \displaystyle \sum_{\beta \in \Phi_{i+1}} x_{\beta}
E_{\beta}$ in $\n{i+1}$ and a root $\alpha$ in $\Phi_i$ with $\alpha \not 
> 2 \alpha_{i+1}$, all in type $D_n$.

\begin{enumerate}
\item \label{basic expansion}
The coefficient of $E_{\alpha}$ in $\Ad \exp X(N)$ is
\[n_{\alpha} 
+ \displaystyle \sum_{\beta \in \Phi^+: \\
\alpha - \beta \in \Phi_{i+1}} m_{\alpha - \beta, \beta} 
x_{\alpha - \beta} n_{\beta}.\]

\item \label{D_n coefficients part 2} If $x_{\beta}$ is zero for each
$\beta < \alpha$ then the coefficient of $E_{\alpha}$ in $\Ad \exp X(N)$
is $n_{\alpha}$.

\end{enumerate}
\end{lemma}

\begin{proof}
By Proposition \ref{near-linearity} the projection 
\[\rho_i \Ad \exp X(N) = 
\rho_i N + \rho_i [X,N] + \rho_i \frac{1}{2}[X,[X,N]].\]  
The coefficient of $E_{\alpha}$ in this expansion is 
\[\displaystyle \begin{array}{ll} n_{\alpha} &+
\displaystyle \sum_{\beta \in \Phi_{i+1}: \alpha - \beta \in \Phi^+} 
m_{\beta, \alpha - \beta} x_{\beta} n_{\alpha -\beta}
\\ &+ 
\displaystyle \sum_{\beta_1, \beta_2 \in \Phi_{i+1}:
\alpha - \beta_1 - \beta_2 \in \Phi^+, \alpha - \beta_1 \in \Phi^+}
\frac{c_{\alpha, \beta_1, \beta_2}
x_{\beta_1} x_{\beta_2} n_{\alpha - \beta_1 - \beta_2}}{2}
\end{array}\]
for  nonzero $c_{\alpha, \beta_1, \beta_2}$ 
determined by the $m_{ \beta, \alpha - \beta} $.  
Since $\alpha$ is not greater than $2 \alpha_{i+1}$ 
the difference $\alpha - \beta_1 - \beta_2$ is not 
positive for any $\beta_1$, $\beta_2$ in $\Phi_{i+1}$.  
Thus the projection $\rho_i(\Ad \exp X(N))$ simplifies to 
$\rho_i (N + [X,N])$, expanded in Part \ref{basic expansion}.
Part \ref{D_n coefficients part 2} follows immediately.
\end{proof}

In the next lemma,  retain the assumption 
that the basis vectors $\{E_{\alpha}: \alpha \in \Phi_i\}$ in the 
$i^{th}$ row are ordered by height from highest to lowest, with an 
arbitrary order fixed for the two roots of same height in type $D_n$.

\begin{lemma} \label{containment}
Fix a regular nilpotent element $N$ in $\n{}$, 
a Hessenberg space $H$, and a Weyl group element $w$
so that $\Ad w^{-1} (E_{\alpha_j}) \in H$ for each simple root $\alpha_j$.
If $\alpha$ is in $\Phi_i$ and $E_{\alpha}$ is not in 
$\Ad w (H)$ then the first nonzero entry in the $\alpha$ row of
$\psi_i(N)$ is $m_{ \alpha_j, \alpha - \alpha_j}n_{\alpha_j}$ 
for some simple root $\alpha_j$.  Furthermore
if $E_{\alpha - \alpha_j}$ is any basis vector whose root differs from
$\alpha$ by a simple root then $E_{\alpha - \alpha_j}$ is 
in $\n{w} \cap \n{i}$.
\end{lemma}

\begin{proof}
Lemma \ref{matrix entries} shows that 
the first entry in the $\alpha$ row that can be nonzero occurs in 
the columns $\beta$ for which $\alpha - \beta$ is as small a positive
root as possible, namely when $\alpha - \beta$ is simple.  
The root $\alpha$ cannot be $\alpha_i$ because $E_{\alpha}$ 
is not in $\Ad w (H)$.  At least one root
$\beta < \alpha$ in $\Phi_i$ differs from $\alpha$ by a simple root
$\alpha_j$, by inspection of Table \ref{classical rows}.  
The entry of $\psi_i(N)$ is $m_{ \alpha_j, \alpha - \alpha_j}n_{\alpha_j}$, 
which is nonzero because $N$ is regular, by Lemma \ref{reg nilps}.

We now show that $\alpha - \alpha_j$ is in $\Phi_w$ for any such
$\alpha_j$.  By hypothesis
$E_{w^{-1} \alpha_j}$ is in $H$ but $E_{w^{-1}(\alpha - \alpha_j) +
w^{-1}\alpha_j}$ is not.  Since $H$ is closed under bracket with
$\mathfrak{b}$, the root space $E_{w^{-1}(\alpha - \alpha_j)}$ is 
in the opposite Borel $\mathfrak{b}^-$ and so $w^{-1}(\alpha - \alpha_j)$ 
is negative.
\end{proof}

\section{Iterated towers of affine fiber bundles}

We are now ready to prove the main lemmata of the paper.  They 
construct affine spaces which are pieces of the intersection of
a Schubert cell with $\H(N,H)$.  The main theorem then uses these
affine spaces to show that each Schubert cell in $\H(N,H)$
has the structure of an iterated tower of affine fiber bundles.

Each Hessenberg space $H$ is the direct
sum of root spaces \cite[Lemma 1]{dMPS}.  We define $\Phi_H$ to be
the set of roots such that $H = \mathfrak{t} \oplus \bigoplus_{\alpha \in \Phi_H} \g{\alpha}$.
Using this correspondance, define $H^c = 
\bigoplus_{\alpha \in \Phi_H^c} \g{\alpha}$ to be the 
complementary sum of root spaces.  
$H^c$ is an ad-nilpotent ideal inside $\mathfrak{b}^-$, 
as discussed in \cite[Section 10]{ST}.

\begin{lemma} \label{main lemma}
Fix $N$ in $\n{}$, $w$ in $W$, and a Hessenberg space $H$ so that
 each simple root vector $E_{\alpha_j}$ is in $\Ad w (H)$.
In types $A_n$, $B_n$, and $C_n$,
the set \[\mathcal{X}_i(N) =
\{X \in \n{i} \cap \n{w}: \rho_i \Ad \exp X (N) \in \rho_i \Ad w (H)\}\]
is homeomorphic to an affine space 
of dimension $|\Phi_w \cap \Phi_i \cap w \Phi_H|$.
\end{lemma}

\begin{proof}
The image $\Ad w(H)$ is the direct sum of root spaces and the Cartan subalgebra, so $\rho_i \Ad w(H) = \Ad w(H) \cap \n{i}$.  By definition of $\theta_i(N)$, the preimage $\theta_i(N)^{-1} (\Ad w (H) \cap \n{i}) = \{X:  \rho_i \Ad \exp X(N) \in  \Ad w(H) \cap \n{i}\}$.  It follows that $\mathcal{X}_i(N) = \n{i} \cap \n{w} \cap \theta_i(N)^{-1} (\Ad w (H) \cap \n{i})$.

In types $A_n$ and $B_n$, the map $\theta_i(N)$ is affine
by Proposition \ref{near-linearity} and so the preimage of the linear
subspace $\Ad w (H) \cap \n{i}$ of $\n{i}$ is affine.
The intersection $\mathcal{X}_i(N)$ of this affine preimage with the linear 
subspace $\n{i} \cap \n{w}$ is also affine.

In type $C_n$, recall that $\gamma_i = \sum_{j=i}^{n-1} 2 \alpha_j +
\alpha_n$ and consider the commutative diagram
\[\begin{array}{rcl}
\n{w} \cap \n{i} & \stackrel{\theta_i(N)}{\longrightarrow} & \n{i} \\
\downarrow & &\downarrow \\
(\n{w} \cap \n{i}) / \g{\gamma_i - \alpha_i} & 
   \stackrel{\stackrel{\sim}{\theta_i}(N)}{\longrightarrow} & 
  \n{i}/\g{\gamma_i}  \end{array}\]
whose vertical arrows are vector space quotients by the root space.  
Define the map $\stackrel{\sim}{\theta_i}(N)$ so the diagram commutes.  It is well-defined because the image 
$\theta_i(N) (\g{\gamma_i-\alpha_i})$ 
is the coset $\rho_i N + \g{\gamma_i} =
\theta_i(N)(0) + \g{\gamma_i}$.
It is affine since the image of $\ad^2 \n{i}$ is 
in $\g{\gamma_i} \oplus \bigoplus_{j<i} \n{j}$, so the preimage
$\stackrel{\sim}{\theta_i}(N)^{-1}(\Ad w (H))$ is affine.

The element $[X]$ in the preimage $\stackrel{\sim}{\theta_i}(N)^{-1} 
(\Ad w (H))$ pulls back to the coset $X + \g{\gamma_i - \alpha_i}$
in $\n{w} \cap \n{i}$.  The image $\theta_i(N)(X + \g{\gamma_i - \alpha_i})$
lies in $\Ad w (H) + \g{\gamma_i}$ by commutativity of the diagram.
Moreover, the restriction of $\theta_i(N)$ to
$X + \g{\gamma_i - \alpha_i}$ is affine because the image $\ad^2 \g{\gamma_i
- \alpha_i}$ is zero.  (The image $\theta_i(N) (X+\g{\gamma_i-\alpha_i})$ intersects $\Ad w (H)$ in exactly one point if $\g{\gamma_i} \not \subseteq \Ad w(H)$ and otherwise is contained in $\Ad w(H)$.)  The fiber over $[X]$ intersects
$\mathcal{X}_i(N)$ in the affine space given
by the preimage of $\theta_i(N)(X+ \g{\gamma_i - \alpha_i}) \cap \Ad w (H)$
under $\theta_i(N)$.  The linear map $\psi_i(N)$ thus determines
the dimension of $\mathcal{X}_i(N)$.

We now show that the dimension of 
each of these fibers in types $A_n$, $B_n$, and $C_n$
is constant if $N$ is regular.  It suffices to study the
linear part of the affine operators and to prove that 
the map $\n{w} \cap \n{i} 
\stackrel{\psi_i(N)}{\longrightarrow} \n{i} \longrightarrow
\Ad w(H^c) \cap \n{i}$ is full rank, since then 
the translation in $\theta_i(N)$ does not affect the dimension of the
kernel.  Since $\n{i} =
(\Ad w (H^c) \oplus \Ad w (H)) \cap \n{i}$ this kernel is 
precisely the preimage of $\Ad w (H) \cap \n{i}$.

The restricted matrix for $\psi_i(N): \n{w} \cap \n{i} \longrightarrow
 \Ad w (H^c) \cap \n{i}$ consists of the rows $\alpha$ in $\Phi_i \cap 
w \Phi_H^c$ and the columns $\beta$ in
$\Phi_w \cap \Phi_i$.  Row $\alpha$ of the full matrix for $\psi_i(N)$
has its first 
nonzero entry $m_{\alpha_j, \alpha - \alpha_j} 
n_{\alpha_j}$ in position $(\alpha, \alpha - \alpha_j)$
for some simple root $\alpha_j$, by Lemma \ref{containment}.  
The columns indexed by
$\alpha - \alpha_j$ are distinct in types $A_n$, $B_n$, and $C_n$ because
the $i^{th}$ row is totally ordered by height in these types.  
Finally, the root $\alpha - \alpha_j$ is in $\Phi_w \cap \Phi_i$
by Lemma \ref{containment}, so the rank of 
the matrix induced by $\psi_i(N)$ on $\n{w} \cap \n{i} \longrightarrow
\Ad w (H^c) \cap \n{i} $ is $|\Phi_i \cap w \Phi_H^c|$.  Note that
$\Phi_i \cap w \Phi_H^c$ is a set of positive roots and
$w^{-1}(\Phi_i \cap w \Phi_H^c)$ is a set of negative roots, since
$H$ contains $\mathfrak{b}$.  This means that
$\Phi_i \cap w \Phi_H^c$ is contained in $\Phi_w \cap \Phi_i$, so 
the kernel of this matrix has dimension 
$|\Phi_w \cap \Phi_i| - |\Phi_i \cap w \Phi_H^c| = | \Phi_w \cap \Phi_i 
\cap w \Phi_H|$.
\end{proof}

The next lemma uses a similar approach for type $D_n$, where
there are  technical difficulties because $\psi_i(N)$ is not
regular nilpotent.  Define a  partition of $\Phi_i$:
\[\begin{array}{c} 
\Phi_i^{0} = \{\alpha \in \Phi_i: \alpha \leq \sum_{j=i}^{n-2} \alpha_j\} \\
\Phi_i^1 = \{\sum_{j=i}^{n-1} \alpha_j, \alpha_n + \sum_{j=i}^{n-2} 
    \alpha_j\}  \\
\Phi_i^{2} = \{\alpha \in \Phi_i: \alpha \geq \sum_{j=i}^{n} 
\alpha_j\}. \end{array}\]
The superscript indicates how many of the simple roots
$\{\alpha_{n-1}, \alpha_n\}$ are summands of the roots in
part $\Phi_i^j$.  Write $\n{i}^j$ for the subspace
$\bigoplus_{\alpha \in \Phi_i^j} \g{\alpha}$, as well as $\rho_i^j$ for the
projection $\n{} \longrightarrow \n{i}^j$ and $\theta_i^j$ 
for the composition $\rho_i^j \circ \theta_i$.  In type $D_n$,
normalize the basis $\{E_{\alpha}\}$ so that
\[\begin{array}{l}m_{\sum_{j=i}^{n-2} \alpha_j, \alpha_{n-1}} 
= m_{\sum_{j=i}^{n-2} \alpha_j, \alpha_{n}} = \\
m_{\alpha_i, \sum_{j=i+1}^{n-1} \alpha_j} =
m_{\alpha_i, \alpha_n+\sum_{j=i+1}^{n-2} \alpha_j} =\\
m_{\sum_{j=i+1}^{n-1} \alpha_j, \alpha_n} =
m_{\alpha_n+\sum_{j=i+1}^{n-2} \alpha_j, \alpha_{n-1}} =1 \end{array}\]
for all $i$ simultaneously.
This is possible by, for instance, \cite[page 54]{Sa}.

This lemma proves that for each $X$ in $\n{i}^0$ and
$Y$ in $\n{i+1}^1 \oplus \n{i+1}^2$, the map
\[\begin{array}{rl}
\theta_D: \n{i}^0 \oplus \n{i+1}^1 \oplus \n{i+1}^2 &\longrightarrow
      \n{i+1}^2 \oplus \n{i}^1 \oplus \n{i}^0 \\
(X,Y) &\mapsto \left( \begin{array}{c}
                \theta_{i+1}^2(N)(Y) \\
                \rho_i^1 \Ad \exp X (\Ad \exp Y(N)) \\
                \theta_i^0(\Ad \exp Y (N))(X) \end{array} \right)
\end{array}\]
is affine and surjects onto the subspace of $\Ad w(H)$ in its image.  
The main step is to 
write the linear part of $\theta_D$ as a matrix
whose first column and last row are zero, and to show that
the remaining minor is block diagonal, each of whose diagonal
blocks is invertible when $N$ is regular.

\begin{lemma} \label{main lemma D}
The set 
\[\mathcal{X}_i(N) = 
\left\{
\begin{array}{l}
(X,Y) \in \n{i}^0 \oplus \n{i+1}^1 
\oplus \n{i+1}^2: \\
\hspace{2em} \theta_D(X,Y) \in
\rho_{i+1}^2\Ad w (H) \oplus \rho_i^1 \Ad w (H) \oplus 
\rho_i^{0} \Ad w (H) 
\end{array} 
\right\}  \]
 is affine of dimension 
\[|\Phi_w \cap (\Phi_i^0 \cup \Phi_{i+1}^1 \cup \Phi_{i+1}^2)| - 
|(\Phi_i^0 \cup \Phi_i^1 \cup \Phi_{i+1}^2) \cap w \Phi_H^c|.\]
\end{lemma}

\begin{proof}
We prove that $\theta_i^0(\Ad \exp Y (N))= \theta_i^0(N)$
for each $Y$ in $\n{i+1}^1 \oplus \n{i+1}^2$ by showing that 
$\Ad \exp Y(N)$ differs from $N$ only in root
spaces which do not affect the map $\theta_i^0$.  Indeed, 
each root in $\Phi_{i+1}^1 \cup \Phi_{i+1}^2$ is greater than at least one of 
$\alpha_{n-1}$ or $\alpha_{n}$, while no root from $\Phi_i^0$ is
greater than either $\alpha_{n-1}$ or $\alpha_n$.
For each $Y$ in $\n{i+1}^1 \oplus \n{i+1}^2$, the coefficient of $E_{\beta}$ 
in the root vector expansion of
$\Ad \exp Y(N)$ agrees with that of $N$ for all $\beta \leq \sum_{j=i}^{n-2}
\alpha_j$, by Lemma \ref{D_n coefficients}.\ref{D_n coefficients part 2}.  
These are the only basis vectors that affect either the translation
or, by Lemma \ref{matrix entries}, the linear part of the affine operators,
 so $\theta_i^0(N) = \theta_i^0(\Ad \exp Y(N))$.

The maps $\theta_i^0(N)$ and $\theta_{i+1}^2(N)$ are
affine operators on  $\n{i}^0$ and $\n{i+1}^1 \oplus \n{i+1}^2$ respectively
by Proposition \ref{near-linearity}.
We write $\rho_i^1 \Ad \exp X (\Ad \exp Y(N))$ explicitly to see that it, too, 
is affine in the $X_{\alpha}$ and $Y_{\alpha}$.
The coefficient of $E_{\sum_{j=i}^{n-1} \alpha_j}$
in $\theta_i (\Ad \exp Y (N))(X)$ is 
\[n_{\sum_{j=i}^{n-1} \alpha_j} + m_{\sum_{j=i+1}^{n-1} \alpha_j, \alpha_i}
n_{\alpha_i} 
Y_{\sum_{j=i+1}^{n-1} \alpha_j} + \sum_{k=i+1}^{n-1} 
m_{\sum_{j=i}^{k-1} \alpha_j, \sum_{j=k}^{n-1}\alpha_j} n_{\sum_{j=k}^{n-1}
\alpha_j} X_{\sum_{j=i}^{k-1} \alpha_j}\] 
by Proposition 
\ref{near-linearity} and Lemma 
\ref{D_n coefficients}.\ref{basic expansion}.  The coefficient of
$E_{\sum_{j=i}^{n-2} \alpha_j + \alpha_n}$ is obtained from this
formula by exchanging 
$\alpha_{n-1}$ and $\alpha_n$.  
Both coefficients are affine functions
in the $X_{\alpha}$ and $Y_{\alpha}$ and so the map
$\theta_D$ is affine.  Since $\mathcal{X}_i(N)$ is
the preimage of the linear space $\Ad w (H) \cap (\n{i}^0 \oplus \n{i}^1
\oplus \n{i+1}^2)$ under this affine map, it is affine itself.

Write the linear part of $\theta_D$ with respect to the
basis of root vectors.  Order the columns
$\Phi_{i+1}^2$ from highest root to 
lowest root, follow with the columns $\Phi_{i+1}^1$ ordered
as in the definition, and then with the
columns $\Phi_i^0$ ordered from highest to lowest.  
Similarly, order
the rows $\Phi_{i+1}^2$, then $\Phi_i^1$, then $\Phi_i^0$, 
within each set ordering by height or by definition. 

The first column corresponds to the root $\alpha_{i+1} + \sum_{j=i+2}^{n-2}
2 \alpha_j + \alpha_{n-1} + \alpha_n$ and the last row to $\alpha_i$.
There is no root in $\Phi_{i+1}^2$,$\Phi_{i+1}^1$,$\Phi_i^1$,
 or $\Phi_i^0$ which is greater
than the former or less than the latter, so this column and row are 
identically zero.

Now examine the minor obtained by omitting the first column and last row.
Form blocks by partitioning the columns into three
sets $\Phi_{i+1}^2$, $\Phi_i^1 \cup \{\sum_{j=i}^{n-2} \alpha_j\}$,
and $\Phi_i^0 \backslash \{\sum_{j=i}^{n-2} \alpha_j\}$, and 
the rows into the sets 
$\Phi_{i+1}^2 \backslash \{\sum_{j=i+1}^n \alpha_j\}$, 
$\Phi_{i+1}^1 \cup \{\sum_{j=i+1}^n \alpha_j\}$, and $\Phi_i^0$.
The matrix is block diagonal because each block to the left of
the block $B$ is indexed  by roots not less than the roots in $B$.
In each diagonal block, 
the first nonzero entry in the row for $\alpha$
is $m_{\alpha - \alpha_j, \alpha_j}n_{\alpha_j}$ and is 
located in a column $\alpha - \alpha_j$ where  $\alpha - \alpha_j$ differs
from $\alpha$ by a simple root.  The root $\alpha - \alpha_j$ is 
unique 
if $\alpha > \sum_{j=i+1}^n \alpha_j$ in $\Phi_{i+1}^2$
or if $\alpha$ is in $\Phi_i^0$
because those root subsets
are totally ordered by height, by inspection of Table \ref{classical rows}.  
Thus, the first and third diagonal blocks are upper-triangular with nonzero
entries along the diagonal.  

Using the previous explicit calculations of $\rho_i^1 \Ad \exp X 
(\Ad \exp Y(N))$  
and computing the coefficient of $E_{\sum_{j=i+1}^n \alpha_j}$ in 
$\theta_{i+1}(N)(Y)$ shows that the second diagonal block is 
\[ \begin{array}{c} \left( \begin{array}{ccc} 
  m_{\sum_{j=i+1}^{n-1} \alpha_j, \alpha_{n}}n_{\alpha_{n}} 
& m_{\alpha_n+\sum_{j=i+1}^{n-2} \alpha_j, \alpha_{n-1}}
n_{\alpha_{n-1}}  & 0\\
m_{\sum_{j=i+1}^{n-1} \alpha_j, \alpha_i}n_{\alpha_i} & 0
& m_{\sum_{j=i}^{n-2} \alpha_j, \alpha_{n-1}}n_{\alpha_{n-1}} \\
0 
& m_{\alpha_n+\sum_{j=i+1}^{n-2} \alpha_j, \alpha_{i}}n_{\alpha_i} 
&  m_{\sum_{j=i}^{n-2} \alpha_j, \alpha_{n}} n_{\alpha_{n}}
\end{array} \right) \\ = \left( \begin{array}{ccc}
   n_{\alpha_n} & n_{\alpha_{n-1}} & 0 \\
   -n_{\alpha_i} & 0 & n_{\alpha_{n-1}} \\
   0 & -n_{\alpha_i} & n_{\alpha_n} \end{array} \right)\end{array} \]
by the basis normalization.
This is invertible since the $n_{\alpha_j}$ are nonzero.  

This confirms that the upper minor of $\theta_D$ is full rank independent
of $X$ and $Y$.
We now prove that the projection of $\theta_D$ to 
$\Ad w (H^c) \cap (\n{i}^0 \oplus \n{i}^1
\oplus \n{i+1}^2)$ is full rank on $\n{w} \cap (\n{i}^0 \oplus
\n{i+1}^1 \oplus \n{i+1}^2)$.  For each row $\alpha$
that is in $w \Phi_H^c$, all of the nonzero columns 
$\alpha - \alpha_j$ for the row are in $\Phi_w$ by Lemma 
\ref{containment}.  This shows that the entries used to
determine the rank of the full matrix are also 
in the matrix restricted and projected to
$\n{w} \cap (\n{i}^0 \oplus \n{i+1}^1 \oplus \n{i+1}^2) \longrightarrow
\Ad w (H^c) \cap (\n{i}^0 \oplus \n{i}^1 \oplus \n{i+1}^2)$.  
Consequently, this
restriction is full rank and so the dimension of its kernel is
$|\Phi_w \cap (\Phi_i^0 \cup \Phi_{i+1}^1 \cup \Phi_{i+1}^2)| - 
|(\Phi_i^0 \cup \Phi_i^1 \cup \Phi_{i+1}^2) \cap w \Phi_H^c|$.
\end{proof}

The main theorem studies the intersection of the
Hessenberg variety $\H(N,H)$ for $N$ in $\mathfrak{b}$
with the Bruhat decomposition of $G/B$ given by the
Borel subgroup corresponding to $\mathfrak{b}$.  The previous
lemmata will
show that each nonempty cell in this decomposition has the
structure of
an iterated affine fiber bundle and so is homeomorphic to an affine
cell.

\begin{theorem} \label{main theorem}
Let $N$ be a regular nilpotent element in $\mathfrak{b}$, $H$ a
Hessenberg space with respect to $\mathfrak{b}$, and 
$\H(N,H)$ the corresponding Hessenberg variety.  
Let $P_w = \H(N,H) \cap BwB/B$ be the intersection of the Schubert cell
corresponding to $w$ with the Hessenberg variety.
The $\{P_w\}$ form a paving by affines
of $\H(N,H)$ when ordered subordinate to the partial order
determined by the length of $w$.  The cell $P_w$ 
is nonempty if and only if $w^{-1} \alpha_i \in
\Phi_H$ for each simple root $\alpha_i$.  If $P_w$ is nonempty
its dimension is $|\Phi_w \cap w \Phi_H|$.
\end{theorem}

\begin{proof}
The $\{P_w\}$ form a paving under any order that respects the length
partial order because the Bruhat decomposition of $G/B$ is a paving
and $\H(N,H)$ is closed in $G/B$.

Consider the set of Lie algebra elements $\Ad u(N)$ for $u \in U_w$.
If $u$ is in $U$ then $\Ad u(N)$ is regular nilpotent
and the coefficient of each simple root vector is nonzero, 
by Lemma \ref{reg nilps}.  Thus,
the element $\Ad u(N)$ can only be in $\Ad w (H)$ if each
$w^{-1} \alpha_i$ is in $\Phi_H$.

We now prove that the condition is sufficient for $P_w$ to be nonempty
and compute the dimension of the
affine cell. Write $U = U_1 \cdots U_n$, factor
$u = u_1 \cdots u_n$ accordingly, and let $X_i$ be an element of $\n{i}$ with 
$\exp X_i = u_i$.  Then $\Ad u(N)$ is 
\[\Ad u(N) = \Ad \exp X_1 (\Ad \exp X_2 (\cdots 
        \Ad \exp X_n (N) \cdots )).\]
(In type $D_n$, further decompose $X_i$ into 
$X_i^0 \in \n{i}^0$ and $X_i^1 \in \n{i}^1 \oplus \n{i}^2$, 
and factor
\[\Ad u(N) = \Ad \exp X_1^1 (\Ad \exp X_1^0( \Ad \exp X_2^1 (\cdots 
        \Ad \exp X_n^0 (N) \cdots ).)\]

Define $\mathcal{Z}_i$ to be the set 
\[\{u_i u_{i+1} \cdots u_n:  u_j \in U_w \cap U_j \hspace{.5em} \forall j,
\rho_j \Ad ( u_i \cdots u_n) (N) \in \rho_j \Ad w (H) \hspace{.5em} 
\forall j \geq i\}.\]
The set $\mathcal{Z}_1$ is homeomorphic to the cell $P_w$
via the map that sends $u_1 \cdots u_n$ to the flag corresponding to 
$(u_1 \cdots u_n)^{-1}w$.
There is a natural map $\mathcal{Z}_i \longrightarrow \mathcal{Z}_{i+1}$
given by $u_i u_{i+1} \cdots u_n \mapsto u_{i+1} \cdots u_n$.  The fiber
over the point $u'$ is the set 
\[\{\exp X_i: X_i \in \n{w} \cap \n{i},
\rho_i \Ad \exp X_i (\Ad u' (N)) \in \rho_i \Ad w (H)\},\] 
namely $\exp \mathcal{X}_i(\Ad u'(N))$ 
of Lemma \ref{main lemma}.  The exponential map
is a homeomorphism on $\n{i}$ and so by Lemma \ref{main lemma} the map
$\mathcal{Z}_i \longrightarrow \mathcal{Z}_{i+1}$ is an affine fiber
bundle of rank $|\Phi_w \cap \Phi_i \cap w \Phi_H|$ for each $i$.
The space $P_w \cong \mathcal{Z}_1$ is thus an iterated tower of
affine fiber bundles and is itself homeomorphic to an affine 
space of dimension $|\Phi_w \cap w \Phi_H|$.

(In type $D_n$ we use the sets
\[\mathcal{Z}_i=\left\{u_i^0 u_{i+1}^1 u_{i+1}^0 \cdots u_n^1:  
\begin{array}{l}u_j^0 \in U_w \cap U_j^0 \hspace{.5em} \forall j, 
  u_j^1 \in U_w \cap (U_j^1 U_j^2) \hspace{.5em} \forall j;\\
\vspace{-.75em}\\
\rho_i^k \Ad (u_i^0 \cdots u_n^1) (N) \in \rho_i^k
\Ad w (H) \textup{ for } k=0,1;\\
\vspace{-.75em}\\
\rho_j \Ad ( u_i^0 \cdots u_n^1) (N) \in \rho_j
\Ad w (H) 
 \hspace{.5em} \forall j > i
\end{array}\right\}.\]
As before, the set $\mathcal{Z}_0$ is homeomorphic to the
cell $P_w$.  The map $\mathcal{Z}_i \longrightarrow
\mathcal{Z}_{i+1}$ given by $(u_i^0 u_{i+1}^1) u'
\mapsto u'$ is an affine fiber bundle
with fiber $\exp \mathcal{X}_i(\Ad  u' (N))$ 
from Lemma \ref{main lemma D}.
This displays $\mathcal{Z}_0$ as an iterated tower of affine
fiber bundles, and so $P_w$ is homeomorphic to affine space.  
The dimension of each fiber is
\[|\Phi_w \cap (\Phi_i^0 \cup \Phi_{i+1}^1 \cup \Phi_{i+1}^2)| - 
|(\Phi_i^0 \cup \Phi_i^1 \cup \Phi_{i+1}^2) \cap w \Phi_H^c|.\]
Summing over $i$ gives $|\Phi_w| - |\Phi^+ \cap w \Phi_H^c|$.  Since
each root in $\Phi^+ \cap w \Phi_H^c$ is in $\Phi_w$, the total dimension is
$|\Phi_w \cap w \Phi_H|$.)
\end{proof}

The statement of this theorem 
is more concise when $N$ is the sum of simple root vectors.

\begin{corollary}
Let $N$ be the sum of simple root vectors $N = \sum_{\alpha_i}E_{\alpha_i}$
in $\mathfrak{b}$, $H$ a
Hessenberg space with respect to $\mathfrak{b}$, and 
$\H(N,H)$ the corresponding Hessenberg variety.  Let $P_w = \H(N,H)
\cap BwB/B$ be the intersection of the Schubert cell corresponding to 
$w$ with the Hessenberg variety.  The $\{P_w\}$ form a paving by
affines of $\H(N,H)$ when ordered subordinate to the partial order
determined by the length of $w$.
The cell $P_w$ is nonempty if and only if $\Ad w^{-1}(N) \in
H$.  If $P_w$ is nonempty
its dimension is $|\Phi_w \cap w \Phi_H|$.
\end{corollary}

\begin{proof}
The sum of simple root vectors 
is regular in all classical types by Lemma \ref{reg nilps}.  Since 
$\Ad w^{-1}(N) = \sum_{\alpha_i} E_{w^{-1} \alpha_i}$ and since
$H$ is a sum of root spaces, the condition
$\Ad w^{-1}(N) \in H$ is equivalent to $w^{-1} \alpha_i \in
\Phi_H$ for each simple $\alpha_i$.
\end{proof}

While the criterion for nonemptiness is more complicated, this
theorem also proves that all regular nilpotent Hessenberg varieties
are paved by affine cells.

\begin{corollary} \label{generic paving}
Fix $\g{}$ of classical type,
let $N$ be a regular nilpotent in $\g{}$, and let $H$ be
a Hessenberg space with respect to $\mathfrak{b}$.  The
Hessenberg variety $\H(N,H)$ is paved by affines.
\end{corollary}

\begin{proof}
Choose an element $\Ad g^{-1}(N)$
in the regular nilpotent orbit which is also in $\n{}$.  
The variety $\H(\Ad g^{-1}(N), H)$ is paved by affines $\{P_w\}$ by
Theorem \ref{main theorem}.
Note that $\H(\Ad g^{-1}(N), H) = g^{-1} \H(N,H)$ and that translation
is a homeomorphism in $G/B$.  This means $\H(N,H)$ is paved by the affine
cells $g P_w$.
\end{proof}

The existence of a paving by affines shows the following.

\begin{corollary} \label{zero odd}
In classical types, for any Hessenberg space $H$ with respect to
$\mathfrak{b}$, the regular nilpotent Hessenberg variety
$\H(N,H)$ has no odd-dimensional cohomology.
\end{corollary}

\begin{proof}
The existence of a paving by complex affine cells
means that the odd-dimensional cohomology of $\H(N,H)$ vanishes
by Lemma \ref{odd-dim}.
\end{proof}

We remark that since $\H(N,H)$ has no odd-dimensional cohomology,
it is equivariantly formal with respect to any algebraic torus action
\cite[p.26]{GKM}.

\end{document}